\documentclass{article}
\usepackage{graphicx}
\usepackage{amssymb,amsmath}
\usepackage{multirow}
\date{}

\makeatletter \@addtoreset{equation}{section}
\renewcommand{\theequation}{\thesection.\@arabic\c@equation}
\newcommand{\affiliation}[1]{\let\thefootnote\relax\footnote{\mbox{}\\ \noindent {#1}}}
\newcommand\bigzero{\makebox(0,0){\text{\huge0}}}

\makeatother

\begin{document}
\title{\textbf{On the Inverting of A General Heptadiagonal Matrix}}

\author{ A. A. KARAWIA\footnote{ Home Address: Mathematics Department, Faculty of Science, Mansoura University,
Mansoura 35516, Egypt. E-mail:abibka@mans.edu.eg}\\
Computer science unit, Deanship of educational services, Qassim University,\\
 P.O.Box 6595, Buraidah 51452, Saudi Arabia. \\
              E-mail: kraoieh@qu.edu.sa
}

\maketitle
%\large{
\begin{abstract}
In this paper, we developed new numeric and symbolic algorithms to find the inverse of any nonsingular heptadiagonal matrix. Symbolic algorithm will not break and it is without setting any restrictive conditions. The computational cost of our algorithms is $O(n)$. The algorithms are suitable for implementation using computer algebra system such as MAPLE, MATLAB and MATHEMATICA. Examples are given to illustrate the efficiency of the algorithms.\bigskip
\end{abstract}

\begin{flushleft}\footnotesize
\hspace*{0.9cm}{\textbf{Keywords}:Heptadiagonal matrices; LU factorization; Determinants; Computer algebra systems(CAS). \\
\hspace*{0.9cm}{\textbf{AMS Subject Classification}:15A15; 15A23; 68W30; 11Y05; 33F10; F.2.1; G.1.0.\\}
}

\end{flushleft}
\newtheorem{alg}{Algorithm}[section]

\section{Introduction}

\hspace*{0.5cm}The$n\times n$ general heptadiagonal matrices take the form:

\begin{equation}
H=\left(
      \begin{array}{cccccccccc}
        d_1 & e_1 & f_1 & g_1 &  & &  &  &  &  \\
        c_2 & d_2 & e_2 & f_2 & g_2 & & & & & \\
        b_3 & c_3 & d_3 & e_3 & f_3 & g_3 & &\bigzero & & \\
        a_4 & b_4 & c_4 & d_4 & e_4 & f_4 & g_4 & & & \\
        & & & & & & & & & \\
        &\ddots &\ddots &\ddots &\ddots &\ddots &\ddots &\ddots & & \\
        & & & & & & & & & \\
         &  &  & a_{n-3} & b_{n-3} & c_{n-3} & d_{n-3} & e_{n-3} & f_{n-3} & g_{n-3} \\
        &  & \bigzero &  & a_{n-2} & b_{n-2} & c_{n-2} & d_{n-2} & e_{n-2} & f_{n-2} \\
        & & & & & a_{n-1} & b_{n-1} & c_{n-1} & d_{n-1} & e_{n-1} \\
        & & & & &  & a_n & b_n & c_n & d_n
      \end{array}
    \right), n>4.
    \end{equation}
where $\{a_i\}_{4\leq i\leq n}$, $\{b_i\}_{3\leq i\leq n}$, $\{c_i\}_{2\leq i\leq n}$, $\{d_i\}_{1\leq i\leq n}$, $\{e_i\}_{1\leq i\leq n}$, $\{f_i\}_{1\leq i\leq n}$, and $\{g_i\}_{1\leq i\leq n}$ are sequences of numbers such that $g_i\neq 0,\quad g_{n-2}=g_{n-1}=g_n=1$ and $f_{n-1}=f_n=e_n=0$.

Heptadiagonal matrices are frequently arise from boundary value problems. So a good technique for compute the inverse of such matrices is required. Also, These kind of matrices appear in many areas of science and engineering[1-9]. To the best of our knowledge, the inversion of a general heptadiagonal matrix of the form (1.1) has not been considered.\\

In [10], Karawia described a reliable symbolic computational algorithm for inverting general cyclic heptadiagonal matrices by using parallel
computing along with recursion. An explicit formula for the determinant of a heptadiagonal symmetric matrix is given in[6]. Many researchers are studied special cases of heptadiagonal matrix. In [11], the authors presented a symbolic algorithm for finding the inverse of any general nonsingular tridiagonal
matrix. A new efficient computational algorithm to find the inverse of a general tridiagonal matrix is presented in [12] based on the Doolittle LU
factorization. In [13], the authors introduced a computationally efficient algorithm for obtaining the inverse of a tridiagonal matrix and a pentadiogonal matrix and they assumed a few conditions to avoid failure in their own algorithm. The motivation of the current paper is to establish efficient algorithms for inverting heptadiagonal matrix. We generalized the algorithm[13] for finding the inverse of a general heptadiagonal matrix and we presented an efficient symbolic algorithm for finding the inverse of such matrices. The development of a symbolic algorithm is considered in order to remove all cases where the numeric algorithm fails.\\

The paper is organized as follows: In Section 2, Main result is presented. New numeric and symbolic algorithms are given in Section 3. In Section 4, Illustrative examples are presented. Conclusions of the work are given in Section 5.\\

\section{Main Result}

In this section, we present recurrence formulas for the columns of the inverse of a heptadiagonal matrix $H$.\\
When the matrix $H$ is nonsingular, its inversion is computed as follows.\\
Let
$$H^{-1}=[S_{ij}]_{1\leq i,j\leq n}=[C_1, C_2, ..., C_n]$$\\
where $C_k$ is the $k$th column of the inverse matrix $H^{-1}$.\\
By using the fact $HH^{-1}=I_n$, where $I_n$ is the identity matrix, the first $(n-3)$ columns can be obtain by relations\\

\begin{equation}
\left.\begin{aligned}
C_{n-3}&=\frac{1}{g_{n-3}}(E_n-d_nC_n-e_{n-1}C_{n-1}-f_{n-2}C_{n-2}),\\
C_{n-4}&=\frac{1}{g_{n-4}}(E_{n-1}-c_nC_n-d_{n-1}C_{n-1}-e_{n-2}C_{n-2}-f_{n-3}C_{n-3}), \\
C_{n-5}&=\frac{1}{g_{n-5}}(E_{n-2}-b_nC_n-c_{n-1}C_{n-1}-d_{n-2}C_{n-2}-e_{n-3}C_{n-3}-f_{n-4}C_{n-4}), \\
C_j&=\frac{1}{g_j}(E_{j+3}-a_{j+6}C_{j+6}-b_{j+5}C_{j+5}-c_{j+4}C_{j+4}-d_{j+3}C_{j+3}-e_{j+2}C_{j+2}-f_{j+1}C_{j+1}),\\
 & \hspace*{10cm} j=n-6, n-7, ...,1,
\end{aligned}
\right\}
\end{equation}
where $E_k$ is the $k$th unit vector.\\

From (2.1), we note that if we knowing the last three columns $C_n, C_{n-1},$ and $C_{n-2}$ then we can recursively compute the remaining $(n-3)$ columns $C_{n-3}, C_{n-4}, ..., C_1$.\\

At this point it is convenient to give recurrence formulas for computing $C_n, C_{n-1},$ and $C_{n-2}$.\\

Consider the sequence of numbers $\{A_i\}_{1\leq i\leq n+3}$, $\{B_i\}_{1\leq i\leq n+3}$,and $\{C_i\}_{1\leq i\leq n+3}$ characterized by a term recurrence relations

\begin{equation}
\left.\begin{aligned}
&A_1=0,\\
&A_2=0,\\
&A_3=1,\\
&d_1A_1+e_1A_2+f_1A_3+g_1A_4=0,\\
&c_2A_1+d_2A_2+e_2A_3+f_2A_4+g_2A_5=0,\\
&b_3A_1+c_3A_2+d_3A_3+e_3A_4+f_3A_5+g_3A_6=0,\\
&a_iA_{i-3}+b_iA_{i-2}+c_iA_{i-1}+d_iA_i+e_iA_{i+1}+f_iA_{i+2}+g_iA_{i+3}=0,\quad i\geq 4,\\
\end{aligned}
\right\}
\end{equation}

\begin{equation}
\left.\begin{aligned}
&B_1=0,\\
&B_2=1,\\
&B_3=0,\\
&d_1B_1+e_1B_2+f_1B_3+g_1B_4=0,\\
&c_2B_1+d_2B_2+e_2B_3+f_2B_4+g_2B_5=0,\\
&b_3B_1+c_3B_2+d_3B_3+e_3B_4+f_3B_5+g_3B_6=0,\\
&a_iB_{i-3}+b_iB_{i-2}+c_iB_{i-1}+d_iB_i+e_iB_{i+1}+f_iB_{i+2}+g_iB_{i+3}=0,\quad i\geq 4,\\
\end{aligned}
\right\}
\end{equation}
and
\begin{equation}
\left.\begin{aligned}
&C_1=1,\\
&C_2=0,\\
&C_3=0,\\
&d_1C_1+e_1C_2+f_1C_3+g_1C_4=0,\\
&c_2C_1+d_2C_2+e_2C_3+f_2C_4+g_2C_5=0,\\
&b_3C_1+c_3C_2+d_3C_3+e_3C_4+f_3C_5+g_3C_6=0,\\
&a_iC_{i-3}+b_iC_{i-2}+c_iC_{i-1}+d_iC_i+e_iC_{i+1}+f_iC_{i+2}+g_iC_{i+3}=0,\quad i\geq 4.\\
\end{aligned}
\right\}
\end{equation}

Now, we can give matrix forms for term recurrences (2.2), (2.3) and (2.4)
\begin{equation}
\mathbf{HA} = -A_{n+1}E_{n-2}-A_{n+2}E_{n-1}-A_{n+3}E_n,
\end{equation}
\begin{equation}
\mathbf{HB} = -B_{n+1}E_{n-2}-B_{n+2}E_{n-1}-B_{n+3}E_n,
\end{equation}
\begin{equation}
\mathbf{HC} = -C_{n+1}E_{n-2}-C_{n+2}E_{n-1}-C_{n+3}E_n,
\end{equation}
where $\mathbf{A}=[A_1,A_2, ..., A_n]^t$, $\mathbf{B}=[B_1,B_2, ..., B_n]^t$, and $\mathbf{C}=[C_1,C_2, ..., C_n]^t$. \\

Let's define the following determinants:
\begin{equation}
X_i =\left | \begin{array}{ccc}
        A_i & A_{n+2} & A_{n+3} \\
        B_i & B_{n+2} & B_{n+3}\\
        C_i & C_{n+2} & C_{n+3}
      \end{array}\right |, \quad i=1,2, ...,n+1,
\end{equation}

\begin{equation}
Y_i =\left | \begin{array}{ccc}
        A_i & A_{n+1} & A_{n+3} \\
        B_i & B_{n+1} & B_{n+3}\\
        C_i & C_{n+1} & C_{n+3}
      \end{array}\right |, \quad i=1,2, ...,n+2,
\end{equation}

\begin{equation}
Z_i =\left | \begin{array}{ccc}
        A_i & A_{n+1} & A_{n+2} \\
        B_i & B_{n+1} & B_{n+2}\\
        C_i & C_{n+1} & C_{n+2}
      \end{array}\right |, \quad i=1,2, ...,n+3.
\end{equation}
By simple calculations, we have
\begin{equation}
\mathbf{HX} = -X_{n+1}E_{n-2},
\end{equation}

\begin{equation}
\mathbf{HY} = -Y_{n+2}E_{n-1},
\end{equation}

\begin{equation}
\mathbf{HZ} = -Z_{n+3}E_n,
\end{equation}
where $\mathbf{X}=[X_1,X_2, ..., X_n]^t$, $\mathbf{Y}=[Y_1,Y_2, ..., Y_n]^t$, and $\mathbf{Z}=[Z_1,Z_2, ..., Z_n]^t$. \\
\\
\textbf{Remark 2.1.} $X_{n+1}=-Y_{n+2}=Z_{n+3}$.\\
\\
\textbf{Lemma 2.1.}(generalization version of Lemma 3.1 in [13]) If $X_{n+1}=0$, then the matrix $\mathbf{H}$ is singular.\\
\\
\textbf{Proof.} The proof is simple. $\Box$\\
\\
\textbf{Theorem 2.1.}(generalization version of theorem 3.1 in [13]) Suppose that $X_{n+1}\neq 0$, then $\mathbf{H}$ is invertible and
\begin{equation}
C_n =\left [ \frac{-Z_1}{Z_{n+3}}, \frac{-Z_2}{Z_{n+3}}, ..., \frac{-Z_n}{Z_{n+3}}\right ]^t,
\end{equation}

\begin{equation}
C_{n-1} =\left [ \frac{-Y_1}{Y_{n+2}}, \frac{-Y_2}{Y_{n+2}}, ..., \frac{-Y_n}{Y_{n+2}}\right ]^t,
\end{equation}

\begin{equation}
C_{n-2} =\left [ \frac{-X_1}{X_{n+1}}, \frac{-X_2}{X_{n+1}}, ..., \frac{-X_n}{X_{n+1}}\right ]^t.
\end{equation}
\\
\textbf{Proof.} Since $det(\mathbf{H})=-\left (\prod_{i=1}^{n-3}g_i\right ) X_{n+1}\neq 0$, then $\mathbf{H}$ is invertible. From (2.11), (2.12) and (2.13) we obtain $C_n, C_{n-1}$, and $C_{n-2}$. The proof is completed. $\Box$

\section{New numeric and symbolic algorithms for the inverse of heptadiagonal matrix}
\hspace*{0.5cm}In this section, we formulate the result in the previous section . It is a numerical algorithm to compute the inverse of a general heptadiagonal matrix of the form (1.1) when it exists.\\
\noindent\rule{6.4in}{1pt}
\begin{alg}
To find the inverse of heptadiagonal matrix (1.1).\\
\hspace*{2.7cm}let $f_{n-1}=f_n=e_n=0$ and $g_{n-2}=g_{n-1}=g_n$.\\
\noindent\rule{6.4in}{1pt}
\textbf{INPUT:} Order of the matrix $n$ and the components $a_i, b_j, c_k, d_l, e_l, f_l,$ and $g_l$ for $i = 4, 5, . . . , n$,\\ \hspace*{1.7cm}$j = 3, 4, . . . , n$, $k = 2, 3, . . . , n$, and $l = 1, 2, . . . , n$,\\
\textbf{OUTPUT:} The inverse of heptadiagonal matrix $\mathbf{H}^{-1}$.\\
\\
\textbf{Step 1:} Compute the sequence of numbers $A_i, B_i$, and $C_i$ for $i=1, 2, ...,n+3$ using (2.2), (2.3) and\\ \hspace*{1.3cm} (2.4) respectively.\\
\textbf{Step 2:} Compute $X_i,\quad i=1,2, ..., n+1$ using (2.8), $Y_i,\quad i=1,2, ..., n+2$ using (2.9) and $Z_i,$\\ \hspace*{1.3cm} $i=1,2, ..., n+3$ using (2.10).\\
\textbf{Step 3:} Compute the last three columns $C_n, C_{n-1},$ and $C_{n-2}$ using (2.14), (2.15), and (2.16) respectively.\\
\textbf{Step 4:} Compute the remaining (n-3)-columns $C_j, \quad j=n-3, n-4, ...,1$ using (2.1).\\
\textbf{Step 5:} Set $\mathbf{H}^{-1}=[C_1, C_2, ..., C_n].$\\
\noindent\rule{6.4in}{1pt}
\end{alg}

The numeric algorithm 3.1 will be referred to as \textbf{NINVHEPTA} algorithm. The computational cost of \textbf{NINVHEPTA} algorithm is $103n+69$ operations.\\

 As can be easily seen, it breaks down unless the conditions $g_{i}\ne0$ are satisfied for all $i = 1, 2, . . . , n-3$. So the following symbolic algorithm is developed in order to remove the cases where the numeric algorithm fails.\\
\noindent\rule{6.4in}{1pt}
\begin{alg}
To find the inverse of heptadiagonal matrix (1.1).\\
\hspace*{2.7cm}let $f_{n-1}=f_n=e_n=0$ and $g_{n-2}=g_{n-1}=g_n$.\\
\noindent\rule{6.4in}{1pt}
\textbf{INPUT:} Order of the matrix $n$ and the components $a_i, b_j, c_k, d_l, e_l, f_l,$ and $g_l$ for $i = 4, 5, . . . , n$,\\ \hspace*{1.7cm}$j = 3, 4, . . . , n$, $k = 2, 3, . . . , n$, and $l = 1, 2, . . . , n$,\\
\textbf{OUTPUT:} The inverse of heptadiagonal matrix $\mathbf{H}^{-1}$.\\
\\
\textbf{Step 1:} If $g_i=0$ for any $i= 1, 2, . . . , n-3$ set $g_i=t$($t$ is just a symbolic name).\\
\textbf{Step 2:} Compute the sequence of numbers $A_i, B_i$, and $C_i$ for $i=1, 2, ...,n+3$ using (2.2), (2.3) and\\ \hspace*{1.3cm} (2.4) respectively.\\
\textbf{Step 3:} Compute $X_i,\quad i=1,2, ..., n+1$ using (2.8), $Y_i,\quad i=1,2, ..., n+2$ using (2.9) and $Z_i,$\\ \hspace*{1.3cm} $i=1,2, ..., n+3$ using (2.10).\\
\textbf{Step 4:} Compute the last three columns $C_n, C_{n-1},$ and $C_{n-2}$ using (2.14), (2.15), and (2.16) respectively.\\
\textbf{Step 5:} Compute the remaining (n-3)-columns $C_j, \quad j=n-3, n-4, ...,1$ using (2.1).\\
\textbf{Step 6:} Substitute the actual value $t=0$ in all expressions to obtain the elements of columns $C_j$,\\ \hspace*{1.3cm} $j=1,2,...,n$.\\
\textbf{Step 7:} Set $\mathbf{H}^{-1}=[C_1, C_2, ..., C_n].$\\
\noindent\rule{6.4in}{1pt}
\end{alg}

The symbolic algorithm 3.2 will be referred to as \textbf{SINVHEPTA} algorithm. The computational cost of \textbf{SINVHEPTA} algorithm is $103n+69$ operations. Based on \textbf{SINVHEPTA} algorithm, a MAPLE procedure for inverting a general nonsingular heptadiagonal matrix $\textbf{H}$ is listed as an Appendix.\\
\section{ILLUSTRATIVE EXAMPLES}
       In this section we give three examples for the sake of illustration.\\
\\
\textbf{Example 4.1.} (Case I: $g_i\neq 0$ for all $i$)\\ Find the inverse of following $10\times 10$ heptadiagonal matrix

\begin{equation}
H_1=\left[ \begin {array}{cccccccccc} 2&1&4&-1&0&0&0&0&0&0
\\ \noalign{\medskip}5&1&1&2&2&0&0&0&0&0\\ \noalign{\medskip}1&2&-3&2&
7&2&0&0&0&0\\ \noalign{\medskip}6&1&3&2&3&-1&3&0&0&0
\\ \noalign{\medskip}0&1&-1&2&2&-3&4&1&0&0\\ \noalign{\medskip}0&0&4&4
&4&1&2&1&1&0\\ \noalign{\medskip}0&0&0&-1&2&-1&3&-3&2&1
\\ \noalign{\medskip}0&0&0&0&3&1&2&1&11&3\\ \noalign{\medskip}0&0&0&0&0
&4&-3&2&1&1\\ \noalign{\medskip}0&0&0&0&0&0&-7&1&1&2\end {array}
 \right]
\end{equation}
\textbf{Solution:}
By applying the \textbf{NINVHEPTA} algorithm, it yields
\begin{itemize}
    \item \textbf{Step 1:} $A=[0,0,1,4,-\frac{9}{2},{\frac {53}{4}},{\frac {21}{4}},{\frac {83}{4}},-{\frac
                            {93}{2}},{\frac {663}{4}},-{\frac {67}{4}},-198,-269]$,\vspace*{0.3cm}\\
                         \hspace*{1.35cm}$B=[0,1,0,1,-\frac{3}{2},{\frac {13}{4}},{\frac {19}{12}},{\frac {41}{12}},-{
                            \frac {47}{6}},{\frac {341}{12}},-{\frac {53}{12}},-{\frac {107}{3}},-{\frac {124}{3}}]$, and \vspace*{0.3cm}\\
                         \hspace*{1.35cm}$C=[1,0,0,2,-\frac{9}{2},{\frac {53}{4}},{\frac {67}{12}},{\frac {269}{12}},-{
                            \frac {221}{6}},{\frac {1781}{12}},-{\frac {881}{12}},-{\frac {578}{3}},-{\frac {730}{3}}]$.

  \item \textbf{Step 2:} $X=[{\frac {4231}{3}},-{\frac {10942}{3}},-{\frac {2146}{3}},-3688,{\frac {4687}{2}},-{\frac {31741}{12}},-{\frac{19873}{12}},{\frac {51721}{12}},{\frac {19773}{2}},-{\frac {154735}{12}},-{\frac {905413}{12}}]$,\vspace*{0.3cm} \\
      \hspace*{1.3cm}$Y=[{\frac {1983}{4}},-{\frac {62693}{4}},{\frac {11759}{6}},-{\frac {82109}{12}},{\frac {49839}{4}},-{\frac {220819}{12}},-{\frac {141107}{12}},-{\frac {5312}{3}},{\frac {160577}{12}},-{\frac {563539}{12}},0,{\frac {905413}{12}}]$, and\vspace*{0.3cm}\\
      \hspace*{1.43cm}$Z=[{\frac {3325}{12}},-{\frac {33928}{3}},{\frac {21211}{12}},-{\frac {22109}{6}},7763,-{\frac {19327}{2}},-{\frac {84955}{12}},{\frac {50981}{12}},-{\frac {15235}{4}},{\frac {76363}{6}},0,0,-{\frac {905413}{12}}]$.

  \item \textbf{Step 3:} $C_{10}=[{\frac {3325}{905413}},-{\frac {135712}{905413}},{\frac {21211}{905413}},-{\frac {44218}{905413}},{\frac {93156}{905413}},
-{\frac {115962}{905413}},-{\frac {84955}{905413}},{\frac {50981}{905413}},-{\frac {45705}{905413}},{\frac {152726}{905413}}]^t$,\vspace*{0.3cm} \\
      \hspace*{1.5cm}$C_9=[-{\frac {5949}{905413}},{\frac {188079}{905413}},-{\frac {23518}{905413}},{\frac {82109}{905413}},-{\frac {149517}{905413}},
{\frac {220819}{905413}},{\frac {141107}{905413}},{\frac {21248}{905413}},-{\frac {160577}{905413}},{\frac {563539}{905413}}]^t$,and\vspace*{0.3cm}\\
\hspace*{1.5cm}$C_8=[{\frac {16924}{905413}},-{\frac {43768}{905413}},-{\frac {8584}{905413}},-{\frac {44256}{905413}},{\frac {28122}{905413}},-{\frac {31741}{905413}},
-{\frac {19873}{905413}},{\frac {51721}{905413}},{\frac {118638}{905413}},-{\frac {154735}{905413}}]^t$.

\item \textbf{Step 4:} $C_7=[-{\frac {51473}{905413}},{\frac {214649}{905413}},{\frac {6848}{905413}},{\frac {139095}{905413}},-{\frac {121161}{905413}},
{\frac {106328}{905413}},{\frac {88422}{905413}},-{\frac {278373}{905413}},-{\frac {103927}{905413}},{\frac {500627}{905413}},]^t$,\vspace*{0.3cm} \\
      \hspace*{1.3cm}$C_6=[-{\frac {80594}{905413}},-{\frac {217}{905413}},{\frac {83035}{905413}},{\frac {170735}{905413}},-{\frac {10659}{905413}},
{\frac {31638}{905413}},-{\frac {14393}{905413}},-{\frac {84414}{905413}},{\frac {14531}{905413}},-{\frac {15434}{905413}}]^t$,\vspace*{0.3cm} \\
      \hspace*{1.3cm}$C_5=[-{\frac {82176}{905413}},{\frac {447486}{905413}},-{\frac {28082}{905413}},{\frac {170806}{905413}},-{\frac {175068}{905413}},
-{\frac {6589}{905413}},{\frac {102273}{905413}},{\frac {9510}{905413}},-{\frac {78091}{905413}},{\frac {392246}{905413}}]^t$,\vspace*{0.3cm} \\
      \hspace*{1.3cm}$C_4=[{\frac {205297}{905413}},-{\frac {910556}{905413}},{\frac {6935}{905413}},-{\frac {472222}{905413}},{\frac {410790}{905413}},-{\frac {147233}{905413}},{\frac {42741}{905413}},{\frac {427692}{905413}},-{\frac {147953}{905413}},{\frac {9724}{905413}}]^t$,\vspace*{0.3cm} \\
      \hspace*{1.3cm}$C_3=[-{\frac {2619}{905413}},-{\frac {30890}{905413}},-{\frac {25421}{905413}},-{\frac {137812}{905413}},{\frac {172515}{905413}},
-{\frac {19216}{905413}},-{\frac {46090}{905413}},{\frac {62775}{905413}},{\frac {11493}{905413}},-{\frac {198449}{905413}}]^t$ ,\vspace*{0.3cm} \\
      \hspace*{1.3cm}$C_2=[-{\frac {29328}{905413}},{\frac {877900}{905413}},-{\frac {53389}{905413}},{\frac {605688}{905413}},-{\frac {491917}{905413}},
{\frac {172702}{905413}},-{\frac {34896}{905413}},-{\frac {501141}{905413}},{\frac {156703}{905413}},{\frac {50083}{905413}},]^t$, and \vspace*{0.3cm} \\
      \hspace*{1.3cm}$C_1=[-{\frac {88555}{905413}},{\frac {552363}{905413}},{\frac {125378}{905413}},-{\frac {28648}{905413}},-{\frac {88835}{905413}},
{\frac {19552}{905413}},-{\frac {17938}{905413}},-{\frac {61611}{905413}},{\frac {46355}{905413}},-{\frac {55155}{905413}}]^t$.
\item \textbf{Step 5:} $H_1^{-1}=\\\left[ \begin {array}{cccccccccc} -{\frac {88555}{905413}}&-{\frac {
29328}{905413}}&-{\frac {2619}{905413}}&{\frac {205297}{905413}}&-{
\frac {82176}{905413}}&-{\frac {80594}{905413}}&-{\frac {51473}{905413
}}&{\frac {16924}{905413}}&-{\frac {5949}{905413}}&{\frac {3325}{
905413}}\\ \noalign{\medskip}{\frac {552363}{905413}}&{\frac {877900}{
905413}}&-{\frac {30890}{905413}}&-{\frac {910556}{905413}}&{\frac {
447486}{905413}}&-{\frac {217}{905413}}&{\frac {214649}{905413}}&-{
\frac {43768}{905413}}&{\frac {188079}{905413}}&-{\frac {135712}{
905413}}\\ \noalign{\medskip}{\frac {125378}{905413}}&-{\frac {53389}{
905413}}&-{\frac {25421}{905413}}&{\frac {6935}{905413}}&-{\frac {
28082}{905413}}&{\frac {83035}{905413}}&{\frac {6848}{905413}}&-{
\frac {8584}{905413}}&-{\frac {23518}{905413}}&{\frac {21211}{905413}}
\\ \noalign{\medskip}-{\frac {28648}{905413}}&{\frac {605688}{905413}}
&-{\frac {137812}{905413}}&-{\frac {472222}{905413}}&{\frac {170806}{
905413}}&{\frac {170735}{905413}}&{\frac {139095}{905413}}&-{\frac {
44256}{905413}}&{\frac {82109}{905413}}&-{\frac {44218}{905413}}
\\ \noalign{\medskip}-{\frac {88835}{905413}}&-{\frac {491917}{905413}
}&{\frac {172515}{905413}}&{\frac {410790}{905413}}&-{\frac {175068}{
905413}}&-{\frac {10659}{905413}}&-{\frac {121161}{905413}}&{\frac {
28122}{905413}}&-{\frac {149517}{905413}}&{\frac {93156}{905413}}
\\ \noalign{\medskip}{\frac {19552}{905413}}&{\frac {172702}{905413}}&
-{\frac {19216}{905413}}&-{\frac {147233}{905413}}&-{\frac {6589}{
905413}}&{\frac {31638}{905413}}&{\frac {106328}{905413}}&-{\frac {
31741}{905413}}&{\frac {220819}{905413}}&-{\frac {115962}{905413}}
\\ \noalign{\medskip}-{\frac {17938}{905413}}&-{\frac {34896}{905413}}
&-{\frac {46090}{905413}}&{\frac {42741}{905413}}&{\frac {102273}{
905413}}&-{\frac {14393}{905413}}&{\frac {88422}{905413}}&-{\frac {
19873}{905413}}&{\frac {141107}{905413}}&-{\frac {84955}{905413}}
\\ \noalign{\medskip}-{\frac {61611}{905413}}&-{\frac {501141}{905413}
}&{\frac {62775}{905413}}&{\frac {427692}{905413}}&{\frac {9510}{
905413}}&-{\frac {84414}{905413}}&-{\frac {278373}{905413}}&{\frac {
51721}{905413}}&{\frac {21248}{905413}}&{\frac {50981}{905413}}
\\ \noalign{\medskip}{\frac {46355}{905413}}&{\frac {156703}{905413}}&
{\frac {11493}{905413}}&-{\frac {147953}{905413}}&-{\frac {78091}{
905413}}&{\frac {14531}{905413}}&-{\frac {103927}{905413}}&{\frac {
118638}{905413}}&-{\frac {160577}{905413}}&-{\frac {45705}{905413}}
\\ \noalign{\medskip}-{\frac {55155}{905413}}&{\frac {50083}{905413}}&
-{\frac {198449}{905413}}&{\frac {9724}{905413}}&{\frac {392246}{
905413}}&-{\frac {15434}{905413}}&{\frac {500627}{905413}}&-{\frac {
154735}{905413}}&{\frac {563539}{905413}}&{\frac {152726}{905413}}
\end {array} \right]$

\end{itemize}\bigskip
\textbf{Example 4.2.}(Case II: $g_i= 0$ for at least one of $i$)\\Find the inverse of following $5\times 5$ heptadiagonal matrix
\begin{equation}
H_2=\left[ \begin {array}{ccccc} 2&3&4&1&0
\\ \noalign{\medskip}-1&1&-2&3&0
\\ \noalign{\medskip}3&5&1&-1&2
\\ \noalign{\medskip}4&-1&3&2&6
\\ \noalign{\medskip}0&2&1&4&-3\end {array}
 \right]
\end{equation}
\\
\textbf{Solution:}\\
i- By applying the \textbf{NINVHEPTA} algorithm, it breaks down since $g_2=0$.\\
ii- By applying the \textbf{SINVHEPTA} algorithm, it yields
\begin{itemize}
  \item \textbf{Step 1:} $A=[0,0,1,-4,14\,{x}^{-1},-{\frac {5\,x+28}{x}},{\frac {5\,x-84}{x}},3\,{\frac {5\,x+14}{x}}]$,\vspace*{0.3cm}\\
                         \hspace*{1.35cm}$B=[0,1,0,-3,8\,{x}^{-1},-8\,{\frac {x+2}{x}},{\frac {7\,x-48}{x}},2\,{\frac {5\,x+12}{x}}]$, and \vspace*{0.3cm}\\
                         \hspace*{1.35cm}$C=[1,0,0,-2,7\,{x}^{-1},-{\frac {5\,x+14}{x}},-42\,{x}^{-1},{\frac {8\,x+21}{x}}]$.

  \item \textbf{Step 2:} $X=[{\frac {55\,x+294}{x}},{\frac {63+40\,x}{x}},-{\frac {183+56\,x}{x}},-3\,{\frac {2\,x+15}{x}},-79\,{x}^{-1},-{\frac {315\,x+901}{x}}]$,\vspace*{0.3cm} \\
      \hspace*{1.3cm} $Y=[-2\,{\frac {35\,x+88}{x}},7\,{\frac {5\,x+13}{x}},2\,{\frac {7\,x+18}{x}},-{\frac {21\,x+65}{x}},-14\,{x}^{-1},0,{\frac {315\,x+901}{x}}]$, and \vspace*{0.3cm} \\
      \hspace*{1.3cm} $Z=[-{\frac {5\,x-548}{x}},5\,{\frac {5\,x-28}{x}},-{\frac {35\,x+194}{x}},25\,{\frac {3\,x+4}{x}},-325\,{x}^{-1},0,0,-{\frac {315\,x+901}{x}}]$.

  \item \textbf{Step 3:} $C_5=[-{\frac {5\,x-548}{315\,x+901}},5\,{\frac {5\,x-28}{315\,x+901}},-{\frac {35\,x+194}{315\,x+901}},25\,{\frac {3\,x+4}{315\,x+901}},-325\, \frac{ 315}{x+901}]^t$,\vspace*{0.3cm} \\
      \hspace*{1.5cm}$C_4=[2\,{\frac {35\,x+88}{315\,x+901}},-7\,{\frac {5\,x+13}{315\,x+901}},-2\,{\frac {7\,x+18}{315\,x+901}},{\frac {21\,x+65}{315\,x+901}},14\, \frac{315}{x+901}]^t$, and \vspace*{0.3cm} \\
      \hspace*{1.5cm}$C_3=[{\frac {55\,x+294}{315\,x+901}},{\frac {63+40\,x}{315\,x+901}},-{\frac {183+56\,x}{315\,x+901}},-3\,{\frac {2\,x+15}{315\,x+901}}
-79\, \frac{315}{x+901}]^t$.
\item \textbf{Step 4:} $C_2=[-545\, \frac{315}{x+901},205\, \frac{315}{x+901},91\, \frac{315}{x+901},111\, \frac{315}{x+901},315\, \frac{315}{x+901}]^t$, and \vspace*{0.3cm} \\
      \hspace*{1.5cm}$C_1=[-5\,{\frac {13\,x+123}{315\,x+901}},10\,{\frac {x+19}{315\,x+901}},56\,{\frac {2\,x+7}{315\,x+901}},-{\frac {33\,x+7}{315\,x+901}},248\, \frac{315}{x+901}]^t$.
  \item \textbf{Step 5:} $H_2^{-1}=\\\left[ \begin {array}{ccccc} -5\,{\frac {13\,x+123}{315\,x+901}}&-545
\, \frac{315}{x+901}&{\frac {55\,x+294}{315\,x+901}}&2\,{\frac {35\,x+88}{315\,x+901}}&-{\frac {5\,x-548}{315\,x+901}}
\\ \noalign{\medskip}10\,{\frac {x+19}{315\,x+901}}&205\, \frac{315}{x+901}&{\frac {63+40\,x}{315\,x+901}}&-7\,{\frac {5\,x+13
}{315\,x+901}}&5\,{\frac {5\,x-28}{315\,x+901}}\\ \noalign{\medskip}56\,{\frac {2\,x+7}{315\,x+901}}&91\, \frac{315}{x+901}&-{
\frac {183+56\,x}{315\,x+901}}&-2\,{\frac {7\,x+18}{315\,x+901}}&-{\frac {35\,x+194}{315\,x+901}}\\ \noalign{\medskip}-{\frac {33\,x+7}{
315\,x+901}}&111\, \frac{315}{x+901}&-3\,{\frac {2\,x+15}{315\,x+901}}&{\frac {21\,x+65}{315\,x+901}}&25\,{\frac {3\,x+4}{315
\,x+901}}\\ \noalign{\medskip}248\, \frac{315}{x+901}&315\, \frac{315}{x+901}&-79\, \frac{315}{x+901}&14\, \frac{315}{x+901}&-325\, \frac{315}{x+901}\end {array} \right]_{x=0}$

  \item \textbf{Step 6:} $H_2^{-1}=\left[ \begin {array}{ccccc} -{\frac {615}{901}}&-{\frac {545}{901}}&
{\frac {294}{901}}&{\frac {176}{901}}&{\frac {548}{901}}
\\ \noalign{\medskip}{\frac {190}{901}}&{\frac {205}{901}}&{\frac {63}
{901}}&-{\frac {91}{901}}&-{\frac {140}{901}}\\ \noalign{\medskip}{
\frac {392}{901}}&{\frac {91}{901}}&-{\frac {183}{901}}&-{\frac {36}{
901}}&-{\frac {194}{901}}\\ \noalign{\medskip}-{\frac {7}{901}}&{
\frac {111}{901}}&-{\frac {45}{901}}&{\frac {65}{901}}&{\frac {100}{
901}}\\ \noalign{\medskip}{\frac {248}{901}}&{\frac {315}{901}}&-{
\frac {79}{901}}&{\frac {14}{901}}&-{\frac {325}{901}}\end {array}
 \right]$.

\end{itemize}
\textbf{Example 4.3.} We consider the following $n\times n$ heptadiagonal matrix in order to demonstrate the efficiency of \textbf{SINVHEPTA} algorithm.\\
\begin{equation}
H=\left(
      \begin{array}{cccccccccc}
        -2 & -1 & 2 & 1 &  & &  &  &  &  \\
        3 & -2 & -1 & 2 & 1 & & & & & \\
        1 & 3 & -2 & -1 & 2 & 1 & &\bigzero & & \\
        2 & 1 & 3 & -2 & -1 & 2 & 1 & & & \\
        & & & & & & & & & \\
        &\ddots &\ddots &\ddots &\ddots &\ddots &\ddots &\ddots & & \\
        & & & & & & & & & \\
         &  &  & 2 & 1 & 3 & -2 & -1 & 2 & 1 \\
        &  & \bigzero &  & 2 & 1 & 3 & -2 & -1 & 2 \\
        & & & & & 2 & 1 & 3 & -2 & -1 \\
        & & & & &  & 2 & 1 & 3 & -2
      \end{array}
    \right).
    \end{equation}

In Table 1. we give a comparison of the running time between \textbf{SINVHEPTA} algorithm and \textbf{MatrixInverse} function in Maple 13.0 for different orders. It was tested in an Intel(R) Core(TM) i7-4700MQ CPU@2.40GHz 2.40 GHz.\\
Table1.\\
Running time(in Seconds) of proposed algorithm and  \textbf{MatrixInverse} function in Maple 13.0.
$$
\begin{tabular}{|c|c|c|}
  \hline
  n & \textbf{SINVHEPTA} algorithm & \textbf{MatrixInverse} function in Maple 13.0\\ \hline
  100 & 1.282 & 1.672 \\ \hline
  200 & 3.078 & 11.687 \\ \hline
  300 & 9.719 & 40.625 \\ \hline
  500 & 44.312 & 208.968 \\ \hline
  1000 & 439.578 & 2441.188 \\
  \hline
\end{tabular}
$$

\section{CONCLUSIONS}

In this work new numeric and symbolic algorithms have been developed for finding the inverse of any nonsingular heptadiagonal matrix. The algorithms are reliable, computationally efficient and the symbolic algorithm removes the cases where the numeric algorithms fail.\\

\noindent\rule{6.4in}{1pt}
\textbf{Appendix. A MAPLE procedure for inverting a general nonsingular heptadiagonal matrix}\\
\noindent\rule{6.4in}{1pt}
{\small  $>$restart:with(LinearAlgebra):\\
hepta$\_$inv := proc(n::posint,a::vector,b::vector,c::vector,d::vector,e::vector,f::vector,g::vector)\\
local i,j;\\
global A,B,C,X,Y,Z,S,Hinv;\\
A := array(1 .. n+3): B := array(1 .. n+3): C := array(1 .. n+3): X := array(1 .. n+1):\\ Y := array(1 .. n+2): Z := array(1 .. n+3):S:=array(1..n,1..n,sparse):\\
for i from 1 to n-3 do\\
     \hspace*{1cm}if g[i] = 0 then g[i] := x end if\\
end do:\\
A[1]:=0:A[2]:=0:A[3]:=1:B[1]:=0:B[2]:=1:B[3]:=0:C[1]:=1:C[2]:=0:C[3]:=0:\\
A[4]:=-simplify(f[1]/g[1]):A[5]:=-simplify((e[2]+A[4]*f[2])/g[2]):A[6]:=-simplify((d[3]+e[3]*A[4]+A[5]*f[3])/g[3]):\\
B[4]:=-simplify(e[1]/g[1]):B[5]:=-simplify((d[2]+B[4]*f[2])/g[2]):B[6]:=-simplify((c[3]+e[3]*B[4]+B[5]*f[3])/g[3]):\\
C[4]:=-simplify(d[1]/g[1]):C[5]:=-simplify((c[2]+C[4]*f[2])/g[2]):C[6]:=-simplify((b[3]+e[3]*C[4]+C[5]*f[3])/g[3]):\\
for i from 4 to n do\\
    \hspace*{1cm}A[i+3] := -simplify((a[i]*A[i-3]+b[i]*A[i-2]+c[i]*A[i-1]+d[i]*A[i]+e[i]*A[i+1]+f[i]*A[i+2])/g[i]):\\
    \hspace*{1cm}B[i+3] := -simplify((a[i]*B[i-3]+b[i]*B[i-2]+c[i]*B[i-1]+d[i]*B[i]+e[i]*B[i+1]+f[i]*B[i+2])/g[i]):\\
    \hspace*{1cm}C[i+3] := -simplify((a[i]*C[i-3]+b[i]*C[i-2]+c[i]*C[i-1]+d[i]*C[i]+e[i]*C[i+1]+f[i]*C[i+2])/g[i]):\\
end do:\\
i := 'i':\\
for i from 1 to n+1 do\\
    \hspace*{1cm}X[i]:=simplify(Determinant(Matrix([[A[n+3], A[n+2], A[i]], [B[n+3], B[n+2], B[i]], [C[n+3], C[n+2], C[i]]]))):\\
    \hspace*{1cm}Y[i]:=simplify(Determinant(Matrix([[A[n+3], A[n+1], A[i]], [B[n+3], B[n+1], B[i]], [C[n+3], C[n+1], C[i]]]))):\\
    \hspace*{1cm}Z[i]:=simplify(Determinant(Matrix([[A[n+2], A[n+1], A[i]], [B[n+2], B[n+1], B[i]], [C[n+2], C[n+1], C[i]]]))):\\
end do:\\
Y[n+2]:=simplify(Determinant(Matrix([[A[n+3], A[n+1], A[n+2]], [B[n+3], B[n+1], B[n+2]], [C[n+3], C[n+1], \hspace*{2.5cm}C[n+2]]]))):\\
Z[n+2]:=simplify(Determinant(Matrix([[A[n+2], A[n+1], A[n+2]], [B[n+2], B[n], B[n+2]], [C[n+2], C[n],\\ \hspace*{2.5cm}C[n+2]]]))):\\
Z[n+3]:=simplify(Determinant(Matrix([[A[n+2], A[n+1], A[n+3]], [B[n+2], B[n+1], B[n+3]], [C[n+2], C[n+1],\\ \hspace*{2.5cm}C[n+3]]]))):\\
i := 'i':\\
for i from 1 to n do\\
    \hspace*{1cm}S[i,n]:=-Z[i]/Z[n+3]:\\
    \hspace*{1cm}S[i,n-1]:=-Y[i]/Y[n+2]:\\
    \hspace*{1cm}S[i,n-2]:=-X[i]/X[n+1]:\\
end do:\\
i := 'i':\\
for i to n do\\
\hspace*{1cm}S[i, n-3] := -simplify((d[n]*S[i, n]+e[n-1]*S[i, n-1]+f[n-2]*S[i, n-2])/g[n-3]);\\
\hspace*{1cm}if i = n then\\
     \hspace*{2cm}S[i, n-3] := simplify(1/g[n-3]+S[i, n-3]):\\
\hspace*{1cm}end if\\
end do:\\
if n=5 then\\
\hspace*{2cm}    i := 'i':\\
   \hspace*{2cm} for i to n do\\
        \hspace*{2.5cm}S[i, n-4] := -simplify((c[n]*S[i,n]+d[n-1]*S[i, n-1]+e[n-2]*S[i, n-2]+\\ \hspace*{5.5cm}f[n-3]*S[i, n-3])/g[n-4]);\\
        \hspace*{2.5cm}if i = n-1 then\\
           \hspace*{3.5cm}S[i, n-4] := simplify(1/g[n-4]+S[i, n-4])\\
        \hspace*{2.5cm}end if\\
    \hspace*{2cm}end do\\
elif n=6 then\\
 \hspace*{2cm}  i := 'i':\\
  \hspace*{2cm}  for i to n do\\
       \hspace*{2.5cm}S[i, n-5] := -simplify((b[n]*S[i,n]+c[n-1]*S[i,n-1]+d[n-2]*S[i, n-2]+e[n-3]*S[i, n-3]+\\ \hspace*{5.5cm}f[n-3]f[n-4]*S[i, n-4])/g[n-5]);  \\
       \hspace*{2.5cm}if i = n-2 then\\
          \hspace*{3.5cm}S[i, n-5] := simplify(1/g[n-5]+S[i, n-5]):\\
       \hspace*{2.5cm}end if\\
   \hspace*{2cm} end do:\\
else\\
  \hspace*{2cm}  i:='i':\\
  \hspace*{2cm}  for j from n-6 by -1 to 1 do\\
       \hspace*{2.5cm} for i to n do\\
            \hspace*{3cm}S[i, j] := -simplify((a[j+6]*S[i,j+6]+b[j+5]*S[i,j+5]+c[j+4]*S[i,j+4]+d[j+3]*S[i,j+3]+\\ \hspace*{5.5cm}f[n-3]e[j+2]*S[i,j+2]+f[j+1]*S[i, j+1])/g[j]);\\
           \hspace*{3cm} if i = j+3 then\\
                 \hspace*{3.5cm}S[i,j] := simplify(1 / g[j] + S[i,j])\\
           \hspace*{3cm} fi:\\
        \hspace*{2.5cm} od:\\
    \hspace*{2cm} od:\\
fi:\\
Hinv:=evalm(S):\\
eval(evalm(Hinv),x=0):\\
end proc:}\\
\noindent\rule{6.4in}{1pt}

\end{document}